\documentclass{amsproc}

\usepackage[T1]{fontenc}
\usepackage[utf8]{inputenc}

\usepackage{amsthm}
\usepackage{amssymb}
\usepackage{mathtools}
\usepackage{pdfpages}
\usepackage{changes}
\usepackage{tikz-cd}
\usepackage{hyperref}

\newcommand{\N}{\mathbb{N}}
\newcommand{\Z}{\mathbb{Z}}

\newcommand{\Sym}{\mathbb{S}}

\newcommand{\id}{\mathrm{id}}

\makeatletter

\newtheorem*{thm*}{Theorem}

\newtheorem*{conjecture*}{Conjecture}
\newtheorem{problem}{Problem}

\makeatother

\begin{document}

\begin{abstract}
	Braces were introduced by Rump as a generalization of Jacobson radical
	rings.  It turns out that braces allow us to use ring-theoretic and
	group-theoretic methods for studying involutive solutions to the
	Yang--Baxter equation. If braces are replaced by skew braces, then one can
	use similar methods for studying not necessarily involutive solutions.
	Here we collect problems related to (skew) braces and set-theoretic
	solutions to the Yang--Baxter equation.
\end{abstract}


\title{Problems on skew left braces}

\author{Leandro Vendramin}

\address{IMAS--CONICET and Departamento de Matem\'atica, FCEN, Universidad de Buenos Aires, Pabell\'on~1, Ciudad Universitaria, C1428EGA, Buenos Aires, Argentina}
\email{lvendramin@dm.uba.ar}

\maketitle
\setcounter{tocdepth}{1}

\subsection*{Introduction}

In this paper I propose several problems in the theory of skew braces. I hope
that these problems will help to strengthen the interest in the theory of skew
braces and set-theoretic solutions to the Yang--Baxter equation. I have not
attempted to review the general theory. I do not discuss problems related to
homology of solutions~\cite{MR3558231}, semibraces~\cite{MR3649817} or
trusses~\cite{trusses}.  I concentrate only on skew braces and the Yang--Baxter
equation.

In~\cite[\S9]{MR1183474} Drinfeld wrote that maybe it would be interesting to
study set-theoretical solutions to the equation
\begin{equation}
	\label{eq:YB}
	R_{12}R_{13}R_{23}=R_{23}R_{13}R_{12},
\end{equation}
where $X$ is a set, $R\colon X\times X\to X\times X$ and $R_{ij}$ is the map acting
like $R$ on the two factors $i$ and $j$ and leaving the third factor alone.
Writing $r=\tau\circ R$, where $\tau$ is the map $(x,y)\mapsto (y,x)$,
Equation~\eqref{eq:YB} becomes
\begin{equation}
	\label{eq:braid}
	r_{12}r_{23}r_{12}=r_{23}r_{12}r_{23},
\end{equation}
where $r_{12}=r\times\id$ and $r_{23}=\id\times r$. 

Let us say that a pair $(X,r)$ is a \emph{set-theoretic solution} to the Yang--Baxter
equation (YBE) if $X$ is a non-empty set and $r\colon X\times X\to X\times X$ is a
bijective map satisfying~\eqref{eq:braid}. 
If $(X,r)$ is a solution of the YBE we write
\[
	r(x,y)=(\sigma_x(y),\tau_y(x)).
\]
We say that the solution $(X,r)$ is \emph{involutive} if $r^2=\id_{X\times X}$
and \emph{non-degenerate} if $\sigma_x\in\Sym_X$ and $\tau_x\in\Sym_X$ for all
$x\in X$, where $\Sym_X$ denotes the group of bijective maps $X\to X$. By
convention, our \emph{solutions} will always be non-degenerate solutions. 

One tool to study set-theoretic solutions to the YBE is the theory of skew left braces.
Braces were introduced by Rump in~\cite{MR2278047}. 
For a recent survey on braces we refer to~\cite{cedo}. 
Skew braces appeared later
in~\cite{MR3647970}.  A \emph{skew left brace} is a triple $(A,+,\circ)$, where
$(A,+)$ and $(A,\circ)$ are groups 
(not necessarily abelian) such that
\begin{equation}
	\label{eq:left}
	a\circ (b+c)=a\circ b-a+a\circ c
\end{equation}
for all $a,b,c\in A$. 
To define skew right braces one needs to
replace~\eqref{eq:left} by 
\begin{equation}
	\label{eq:right}
	(a+b)\circ c=a\circ c-c+b\circ c.
\end{equation}
The inverse of an element $a\in A$ with respect to the circle (or
multiplicative) group of $A$ will be denoted by $a'$. Examples of skew left braces can
be found in~\cite{MR3320237,MR3177933,MR3647970,MR3291816}. 

A skew \emph{two-sided} brace is a skew left brace that is also a skew right
brace with respect to the same pair of operations.  

If $\mathcal{X}$ is a class of groups, a skew left brace will be called of type
$\mathcal{X}$ if its additive group belongs to $\mathcal{X}$.  Skew left braces
of abelian type are those braces introduced by Rump in~\cite{MR2278047}.  By
convention, a \emph{left brace} will always be a skew left brace of abelian
type, i.e. with abelian additive group. 

Braces and skew left braces have several interesting connections, see for
example~\cite{MR3291816,MR3763907}. The connection between skew left braces and
non-degenerate solutions of the YBE is given in the following theorem,
see~\cite[Theorem 3.1]{MR3647970}:

\begin{thm*}
	If $A$ is a skew left
	brace, then the map
	\[
		r_A\colon A\times A\to A\times A,\quad
		r_A(a,b)=(-a+a\circ b,(-a+a\circ b)'\circ a\circ b),
	\]
	is a non-degenerate
	set-theoretic solution of the YBE. Moreover, $r_A$ is involutive if and
	only if $A$ is of abelian type.
\end{thm*}

Solutions produced from skew left braces are universal. In a
different language, one finds the following result in
\cite{MR1722951,MR1769723,MR1809284}:  

\begin{thm*} 
	If $(X,r)$ be a solution of the YBE, then there
	exists a unique skew left brace structure over the group
	\[
		G=G(X,r)=\langle X:x\circ y=u\circ v\text{ whenever $r(x,y)=(u,v)$}\rangle
	\]
	such that the diagram
	\[
		\begin{tikzcd}
			X\times X\arrow[r, "r"] \arrow["\iota\times\iota"', d]
			& X\times X \arrow[d, "\iota\times\iota" ] \\
			G\times G \arrow[r, "r_{G}"]
			& G\times G
		\end{tikzcd}
	\]
	commutes, where $\iota\colon X\to G(X,r)$ is the canonical map.  
\end{thm*}

One proves that the pair $(G(X,r),\iota)$ satisfies a universal property,
see~\cite[Proposition 3.9]{MR3647970} and~\cite[Theorem 3.5]{MR3763907}.  The
group $G(X,r)$ is infinite since for example the degree map $X\to\Z$,
$x\mapsto 1$, extends to a group homomorphism $G(X,r)\to\Z$.  Several
properties of $G(X,r)$
are discussed in~\cite{MR2301033}.
If the set $X$ of the solution $(X,r)$ is finite, then one can realize this
solution using a finite skew left brace, see~\cite{Bachiller3,MR3177933}.

\subsection*{Problems}

It was observed by Rump~\cite{MR2132760} that non-degenerate involutive
solutions are in bijective correspondence with non-degenerate cycle sets. A
\emph{cycle set} is a pair $(X,\cdot)$, where $X\times X\to X$, $(x,y)\mapsto
x\cdot y$, is a map such that each $\varphi_x\colon X\to X$, $y\mapsto x\cdot
y$, is bijective, and 
\[
(x\cdot y)\cdot (x\cdot z)=(y\cdot
x)\cdot (y\cdot z)
\]
for all $x,y,z\in X$.

A cycle set $(X,\cdot)$ is said to be \emph{non-degenerate} if the map 
$x\mapsto x\cdot x$ is bijective. Finite cycle sets are non-degenerate,
see~\cite[Theorem~2]{MR2132760}. The correspondence between non-degenerate
cycle sets and involutive solutions is given by 
\[
	r(x,y)=((y*x)\cdot y,y*x)
\]
where $x*y=z$ if and only if $x\cdot z=y$. 
Homomorphisms of cycle sets are defined in the usual manner.  

\begin{problem}
	Construct the free cycle set.	
\end{problem}

It would be interesting to have a nice description of the free cycle set. This
nice description could be used for example to compute homology.

In~\cite{MR1722951}, Etingof, Schedler and Soloviev constructed all
non-degenerate involutive solutions of cardinality $\leq8$, see
Table~\ref{tab:numbers}. To construct finite involutive solutions one needs to
construct finite cycle sets. 

\begin{table}[h]
        \centering
        \caption{Involutive non-degenerate solutions of cardinality $\leq8$.}
        \begin{tabular}{|c|c|c|c|c|c|c|c|c|c|c|}
            \hline
            cardinality & 1 & 2 & 3 & 4 & 5 & 6 & 7 & 8\\
			\hline
            solutions & 1 & 2 & 5 & 23 & 88 & 595 & 3456& 34528\\
            square-free & 1 & 1 & 2 & 5 & 17 & 68 & 336& 2041\\
            indecomposable & 1 & 1 & 1 & 5 & 1 & 10 & 1 & 98\\
            \hline
        \end{tabular}
        \label{tab:numbers}
\end{table}

\begin{problem}
	\label{pro:size9}
	Construct all cycle sets of small cardinality.
\end{problem}

One can try for example with the construction of cycle sets of cardinality nine. Of
course Problem~\ref{pro:size9} refers to the construction of all isomorphism
classes of cycle sets.  It would also be interesting to construct all cycle
sets (or involutive solutions) of small cardinality under particular assumptions. For
example, a solution $(X,r)$ of the YBE is said to be \emph{square free} if
$r(x,x)=(x,x)$ for all $x\in X$, this means that one needs cycle sets
$(X,\cdot)$ such that $x\cdot x=x$ for all $x\in X$.

\begin{problem}
	\label{pro:squarefree}
	Construct all square-free cycle sets 
	of small cardinality. 
\end{problem}

An involutive solution $(X,r)$ is said to be \emph{indecomposable} if the group
$\mathcal{G}(X,r)$ defined as the subgroup of $\Sym_X$ generated by
$\{\sigma_x:x\in X\}$ acts transitively on $X$.  Similarly one defines
indecomposable cycle sets.  The classification of indecomposable solutions with
a prime number of elements appears in~\cite{MR1783929,MR1848966}. 

\begin{problem}
	\label{pro:indecomposable}
	Construct all indecomposable
	cycle sets of small cardinality.
\end{problem}

Maybe Problems~\ref{pro:size9}--\ref{pro:indecomposable} and similar problems
could be studied using constraint satisfaction methods.  

The number of finite
involutive solutions increases rapidly with the number of elements of the
underlying set.  Therefore it makes sense to ask for an estimation:

\begin{problem}
	Estimate the number of cycle sets of cardinality $n$ for $n\to\infty$. 
\end{problem}

Much less is known for non-involutive solutions. Before going into the general
problem of constructing non-involutive solutions, one could start with
injective solutions. A solution $(X,r)$ is said to be \emph{injective} if the canonical
map $X\to G(X,r)$ is injective.  As it was proved in~\cite{MR1722951},
involutive solutions are always injective: 
\[
	\{\text{involutive solutions}\}\subsetneq
	\{\text{injective solutions}\}\subsetneq
	\{\text{solutions}\}.
\]

\begin{problem}
	Construct all injective non-involutive solutions of
	low cardinality.
\end{problem}

\begin{problem}
	\label{pro:all}
	Construct all solutions of small cardinality.
\end{problem}

Problem~\ref{pro:all} is related to the construction of finite biquandles of
small cardinality, see for example~\cite{MR2819176,MR3666513}.  
A different approach to construct all solutions could be based on skew left
braces, see~\cite{Bachiller3};  this method
requires the classification of finite skew left braces.

Simple involutive solutions were defined in~\cite[\S2]{MR3437282}. A surjective
map $p\colon X\to Y$ of involutive solutions is said to be a \emph{covering} if
all the fibers $p^{-1}(y)$ have the same cardinality. A covering $X\to Y$ is
said to be \emph{trivial} if either $|Y|=1$ or $|Y|=|X|$.  An involutive
solution $(X,r)$ is said to be \emph{simple} if $|X|>1$ and any covering $X\to
Y$ is trivial.  

\begin{problem}
	Classify finite simple involutive solutions.
\end{problem}

One could simply ask for examples of small involutive simple solutions: 

\begin{problem}
	Construct finite simple involutive solutions of small cardinality.
\end{problem}

It would be interesting to understand the simplicity of an involutive solution in the
language of left braces: 

\begin{problem}
	Is it possible to read off the simplicity of an involutive solution $(X,r)$ in
	terms of the left brace $G(X,r)$?
\end{problem}

It would be also interesting to develop the theory of non-involutive simple
solutions. Some ideas could be obtained if one reads off the simplicity of a
solution in terms of a skew left brace where it is realized. 

\begin{problem}
	\label{pro:torsion}
	Let $(X,r)$ be a solution to the YBE.  When does the multiplicative group of
	$G(X,r)$ have torsion?
\end{problem}

Problem~\ref{pro:torsion} was posed by Bachiller in~\cite{Bachiller3}.
Gateva--Ivanova and Van den Bergh proved in~\cite{MR1637256} that in the case
of involutive solutions the multiplicative group of $G(X,r)$ is always
torsion-free. It is easy to construct examples of non-involutive solutions
$(X,r)$ such that the multiplicative group of $G(X,r)$ has torsion.

\begin{problem}
	Construct the free (skew) left brace.
\end{problem}

Based on the work of Bachiller~\cite{MR3465351} and Catino and
Rizzo~\cite{MR2486886} on regular subgroups, an algorithm to construct finite
skew left braces was introduced in~\cite{MR3647970}.  With some exceptions,
this algorithm was used to construct small (skew) left braces.  For $n\in\N$
let $s(n)$ be the number of non-isomorphic skew left braces of size $n$ and let
$b(n)$ be the number of non-isomorphic left braces of size $n$.  Table
\ref{tab:braces} shows some values of $s(n)$ and $b(n)$ and several open cases.
Maybe some of these open cases are 
fairly straightforward computational projects.  

\begin{table}[h]
	\caption{The number of non-isomorphic left braces and skew left braces.}
	\begin{tabular}{|c|cccccccccccc|}
		\hline
		$n$ & 1 & 2 & 3 & 4 & 5  & 6 & 7 & 8 & 9 & 10 & 11 & 12\tabularnewline
		$b(n)$ & 1  & 1  & 1  & 4  & 1  & 2  & 1  & 27  & 4  & 2  & 1  & 10\tabularnewline
		$s(n)$ & 1 & 1 & 1 & 4 & 1 & 6 & 1 & 47 & 4 & 6 & 1 & 38 \tabularnewline
		\hline
		$n$ & 13 & 14 & 15 & 16 & 17 & 18 & 19 & 20 & 21 & 22 & 23 & 24\tabularnewline
		$b(n)$ & 1  & 2  & 1  & 357  & 1  & 8  & 1  & 11  & 2  & 2  & 1  & 96\tabularnewline
		$s(n)$ & 1 & 6 & 1 & 1605 & 1 & 49 & 1 & 43 & 8 & 6 & 1 & 855 \tabularnewline
		\hline
		$n$ & 25 & 26 & 27 & 28 & 29 & 30 & 31 & 32 & 33 & 34 & 35 & 36\tabularnewline
		$b(n)$ & 4  & 2  & 37  & 9  & 1  & 4 & 1  & ?  & 1  & 2  & 1  & 46\tabularnewline
		$s(n)$ & 4 & 6 & 101 & 29 & 1 & 36 & 1 & ? &  1 & 6 & 1 & 400\tabularnewline
		\hline
		$n$ & 37 & 38 & 39 & 40 & 41 & 42 &  43 & 44 & 45 & 46 & 47 & 48\tabularnewline
		$b(n)$  & 1  & 2  & 2  & 106  & 1  & 6  & 1  & 9  & 4 & 2  & 1  & 1708\tabularnewline
		$s(n)$ & 1 & 6 & 8 & 944 & 1 & 78 & 1 & 29 & 4 & 6 & 1 & 66209 \tabularnewline
		\hline
		$n$    & 49 & 50 & 51 & 52 & 53 & 54 & 55 & 56 & 57 & 58 & 59 & 60\tabularnewline
		$b(n)$ &  4 & 8 & 1 & 11 & 1 & 80 & 2 & 91 & 2 & 2 & 1 & 28\tabularnewline
		$s(n)$ & 4 & 51 & 1 & 43 & 1 & ? &  12 & 815 & 2 & 6 & 1 & 418\tabularnewline

		\hline
		$n$ & 61 & 62 & 63 & 64 & 65 & 66 & 67 & 68 & 69 & 70 & 71 & 72\tabularnewline
		$b(n)$ & 1 & 2 & 11 & ? & 1 & 4 & 1 & 11 & 1 & 4 & 1 & 489\tabularnewline
		$s(n)$ & 1 & 6 & 11 & ? &  1 & 36 & 1 & 43 & 1 & 36 & 1 & 17790 \tabularnewline

		\hline
		$n$ & 73 & 74 & 75 & 76 & 77 & 78 & 79 & 80 & 81 & 82 & 83 & 84\tabularnewline
		$b(n)$ & 1 & 2 & 5 & 9 & 1 & 6 & 1 & 1985 & ? & 2 & 1 & 34\tabularnewline
		$s(n)$ & 1 & 6 & 14 & 29 & 1 & 78 & 1 & ? &  ? &  6 & 1 & 606\tabularnewline

		\hline
		$n$ & 85 & 86 & 87 & 88 & 89 & 90 & 91 & 92 & 93 & 94 & 95 & 96\tabularnewline
		$b(n)$ & 1 & 2 & 1 & 90 & 1 & 16 & 1 & 9 & 2 & 2 & 1 & ?\tabularnewline
		$s(n)$ & 1 & 6 & 1 & 800 & 1 & 294 & 1 & 29 & 8 & 6 & 1 & ? \tabularnewline

		\hline
		$n$ & 97 & 98 & 99 & 100 & 101 & 102 & 103 & 104 & 105 & 106 & 107 & 108\tabularnewline
		$b(n)$ & 1 & 8 & 4 & 51 & 1 & 4 & 1 & 106 & 2 & 2 & 1 & 494 \tabularnewline
		$s(n)$ & 1 & 53 &  4 & 711 & 1 & 36 & 1 & 944 &  8 &  6 &  1 & ? \tabularnewline

		\hline
		$n$ & 109 & 110 & 111 & 112 & 113 & 114 & 115 & 116 & 117 & 118 & 119 & 120\tabularnewline
		$b(n)$ & 1 & 6 & 2 & 1671 & 1 & 6 & 1 & 11 & 11 & 2 & 1 & 395\tabularnewline
		$s(n)$ & 1 & 94 &  8 &  ? &  1 & 78 &  1 & 43 &  47 &  6 &  1 & ? \tabularnewline
		\hline

		$n$ & 121 & 122 & 123 & 124 & 125 & 126 & 127 & 128 & 129 & 130 & 131 & 132\tabularnewline
		$b(n)$ & 4 & 2 & 1 &  9  & 49 & 36 & 1 & ? & 2 &  4 & 1 & 24\tabularnewline
		$s(n)$ & 4 & 6 & 1 & 29 & 213 & 990 & 1 & ? & 8 & 36 & 1 & 324 \tabularnewline
		\hline

		$n$  & 133 & 134 & 135 & 136 & 137 & 138 & 139 & 140 & 141 & 142 & 143 & 144\tabularnewline
		$b(n)$ & 1 & 2 & 37 & 108 & 1 & 4 & 1 & 27 & 1 & 2 & 1 & ?\tabularnewline
		$s(n)$ & 1 & 6 & 101 & ? & 1 & 36 & 1 & 395 & 1 & 6 & 1 & ? \tabularnewline
		\hline
		
		$n$ & 145 & 146 & 147 & 148 & 149 & 150 & 151 & 152 & 153 & 154 & 155 & 156\tabularnewline
		$b(n)$ & 1 & 2 & 9 & 11 & 1 & 19 & 1 & 90 & 4 & 4 & 2 & 40\tabularnewline
		$s(n)$ & 1 & 6 & ? & 43 & 1 & ? & 1 & ? & 4 & 36 & 12 & 782 \tabularnewline
		\hline

		$n$ & 157 & 158 & 159 & 160 & 161 & 162 & 163 & 164 & 165 & 166 & 167 & 168\tabularnewline
		$b(n)$ & 1 & 2 & 1 & ? & 1 & ? & 1 & 11 & 2 & 2 & 1 & 443\tabularnewline
		$s(n)$ & 1 & 6 & 1 & ? & 1 & ? & 1 & 43 & 12 & 6 & 1 & ? \tabularnewline
		\hline
	\end{tabular}
	\label{tab:braces}
\end{table}

\begin{problem}
	\label{pro:32}
	Construct all left braces of size $32$. 
\end{problem}

Some results related to Problem~\ref{pro:32} are shown in
Table~\ref{tab:braces32}.  According to these results, one needs to construct
left braces of size $32$ with additive group isomorphic to $C_2^{5}=C_2\times
C_2\times C_2\times C_2\times C_2$, where $C_2$ denotes the cyclic
(multiplicative) group with two elements.  Apparently, this particular case
cannot be done in reasonable time with the techniques of~\cite{MR3647970}.  The
number of (skew) left braces of size $64$, $96$ or $128$ seems to be extremely
large and our computational methods are not strong enough to construct them
all. 

\begin{table}[h]
	\centering
	\caption{Some calculations of braces of size $32$.}
	\begin{tabular}{|c|c|c|}
		\hline
		additive group & number & time needed\tabularnewline
		\hline
		$C_{32}$ & 9 & 0.1\tabularnewline
		$C_8\times C_4$ & 1334 & 4 hours\tabularnewline
		$C_{16}\times C_2$ &  120 & 30 minutes\tabularnewline
		$C_4^2\times C_2$ & 13512 & 13 days\tabularnewline
		$C_8\times C_2^2$ & 1547 & 10 hours\tabularnewline
		\hline
	\end{tabular}
	\label{tab:braces32}
\end{table}

\begin{problem}
	\label{pro:p^n}
	Let $p$ be a prime number.  Classify 
	skew left braces of size $p^n$. 
\end{problem}

Problem~\ref{pro:p^n} is important because is the key step in the
classification of skew left braces of nilpotent type. 


Left braces of size $p^2$ and $p^3$, where $p$ is a prime number, were
classified by Bachiller in~\cite{MR3320237}. He proved that $b(p^2)=4$ and
$b(p^3)=6p+19$.  Since groups of order $p^2$ are abelian, it follows that
$s(p^2)=4$.  Skew left braces of size $p^3$ were classified by Nejabati Zenouz
in~\cite{NZ}. He proved that $s(p^3)=6p^2+8p+23$ whenever $p>3$. From
Table~\ref{tab:braces} one gets $s(8)=47$ and $s(27)=101$.

\begin{problem}
	\label{pro:pq}
	Let $p$ and $q$ be different prime numbers. 
	Construct all skew left braces of size $pq$. 
\end{problem}

It should be fairly easy to solve Problem~\ref{pro:pq} using the techniques
of~\cite{MR2030805}.

\begin{problem}
	\label{pro:p^2q}
	Let $p$ and $q$ be different prime numbers. Consruct all skew left braces
	of size $p^2q$. 
\end{problem}

Left braces of size $p^2q$ for prime numbers $p$ and $q$ such that $q>p+1$ were
classified by Dietzel in~\cite{Dietzel}. He proved that 
\[
	b(4q)=\begin{cases}
		9 & \text{if $4\nmid q-1$},\\
		11 & \text{if $4\mid q-1$},
	\end{cases}
\]
if $q>3$, and that 
\[
	b(p^2q)=\begin{cases}
		4 & \text{if $p\nmid q-1$},\\
		p+8 & \text{if $p\mid q-1$ and $p^2\nmid q-1$},\\
		2p+8 & \text{if $p^2\mid q-1$},
	\end{cases}
\]
if $q>p+1>3$. 

\begin{problem}
	\label{pro:estimate_braces}
	Estimate the number of (skew) left braces of size $n$ for $n\to\infty$. 
\end{problem}

The following problem appears in~\cite{cedo}: 

\begin{problem}
	\label{pro:aut}
	Compute automorphism groups of skew left braces of size $p^n$. 
\end{problem}

Automorphism groups of skew left braces of size $p^3$
were computed in~\cite{NZpaper}.

In~\cite[2.28(I)]{MR2095675} Gateva--Ivanova made the following conjecture: For
each finite involutive square-free solution $(X,r)$ there exist $x,y\in X$ such
that $x\ne y$ and $\sigma_x=\sigma_y$.  It was proved by Ced\'o, Jespers and
Okni\'nski in~\cite{MR2652212} that the conjecture is true if the group
$\mathcal{G}(X,r)$ generated by $\{\sigma_x:x\in X\}$ is abelian.  Later in~\cite{MR2885602} and
in~\cite{MR3177933} it was shown using other methods that the result is also valid
if $\mathcal{G}(X,r)$ is infinite abelian. In full generality, the conjecture is
now known to be false. The first counterexample appeared in~\cite{MR3437282};
other counterexamples were later constructed in~\cite{MR3719300,MR3780533}.  It
would be interesting to find essentially new counterexamples. 

\begin{problem}
	\label{pro:GI}
	Construct counterexamples to Gateva--Ivanova conjecture. 
\end{problem}

One could ask, for example, for counterexamples of size nine.  Computer
calculations show that there is only one counterexample to Gateva--Ivanova
conjecture among the 38698 involutive solutions of size $\leq8$.  

It could be enlightening to attack Problem~\ref{pro:GI} using the theory of
braces.  Let $A$ be a brace and $X$ be a subset of $A$ such that the
restriction $r_X=r_A|_{X\times X}$ of $r_A$ to $X\times X$ is a solution to the
YBE. We say that $(A,X)$ is a \emph{Gateva--Ivanova pair} if the solution
$(X,r_X)$ is a counterexample to Gateva--Ivanova conjecture.

\begin{problem}
	\label{pro:pairs}
	Find Gateva--Ivanova pairs. 
\end{problem}

One could start studying Problem~\ref{pro:pairs} by inspecting the database of
small braces. I should remark that counterexamples to Gateva--Ivanova
conjecture might provide new examples of 
Artin-Schelter regular algebras with global dimension
$>3$ with interesting properties to study. Gateva--Ivanova conjecture motivated
a deeper study of the structure of braces and related objects, see for
example~\cite{MR3177933, GI15,MR2278047}.

Gateva--Ivanova conjecture is of course related with the retractability of
involutive solutions, introduced in~\cite{MR1722951}.  For an involutive
solution $(X,r)$, consider the equivalence relation on $X$ given by $x\sim y$ if
and only if $\sigma_x=\sigma_y$. 
Denoting by $Y$ the set of equivalence classes of $X$, the 
map $r$ induces a function $\overline{r}\colon Y\times Y\to Y\times Y$. The 
pair $(Y,\overline{r})$ is a solution to the YBE which will be denoted by 
$\operatorname{Ret}(X, r)$.
One defines inductively 
\[
\operatorname{Ret}^{m+1}(X, r) = \operatorname{Ret}(\operatorname{Ret}^m (X, r))
\]
for all $m\geq1$.  The solution $(X, r)$ is said to be a \emph{multipermutation
solution} if there is a positive integer $m$ such that $\operatorname{Ret}^m(X,
r)$ has only one element. The solution $(X,r)$ is said to be
\emph{irretractable} if $\operatorname{Ret}(X, r) = (X, r)$.

A group $G$ admits a \emph{left ordering} if it admits a total
ordering $<$ such that if $x<y$ then $zx<zy$ for all $x,y,z\in G$. 
In~\cite{MR3815290} it is proved that a finite involutive set-theoretic
solution $(X, r)$ of the YBE is a multipermutation solution if and only if its
structure group $G(X,r)$ admits a left ordering. One of the implications was
proved in ~\cite{MR3572046,MR2189580}.

\begin{problem}
	\label{pro:left}
	Let $(X,r)$ be an injective non-involutive solution of the YBE. When does the 
	multiplicative group of $G(X,r)$ admit a left ordering?
\end{problem}

Problem~\ref{pro:left} appears in~\cite{LV3}. 

A group $G$ satisfies the \emph{unique product property} if for every finite
non-empty subsets $A$ and $B$ of $G$ there is an element $x$ which can be
written uniquely as $x = ab$ with $a\in A$ and $b\in B$. 

\begin{problem}
	\label{pro:UPP}
	Let $(X, r)$ be an involutive irretractable solution of the YBE.
	Can $G(X,r)$ satisfy the unique product property?
\end{problem}

In~\cite[Example 8.2.14]{MR2301033} Jespers and Okni\'nksi proved that the
structure group $G(X,r)$ of a certain involutive irretractable solution $(X,r)$
of cardinality four does not satisfy the unique product property.  They showed that
this structure group contains a subgroup isomorphic to the celebrated Promislow
group~\cite{MR940281}. 

A group $G$ is said to be \emph{diffuse} if for every finite non-empty subset
$A$ of $G$ there exists $a\in A$ such that for all $g\in G\setminus\{1\}$,
either $ga\not\in A$ or $g^{-1}a\not\in A$.  Diffuse groups satisfy the unique
product property. However, the precise relation between diffuseness and unique
products is not clear at the moment.  In~\cite{LV3} it is proved that the
structure group $G(X,r)$ of an involutive solution $(X,r)$ is diffuse if and
only if $(X,r)$ is a multipermutation solution.  Therefore a positive answer to
Problem~\ref{pro:UPP} would give an example of a non-diffuse group with the
unique product property.

\begin{problem}
	\label{pro:MP}
	Study non-involutive multipermutation solutions.
\end{problem}

Some work related to Problem~\ref{pro:MP} can be found in~\cite{LV3,MR3763907}.

As it was observed by Rump, two-sided braces are equivalent to radical rings.
The multiplication of this radical ring is the operation 
\[
a*b=-a+a\circ b-b.
\]
It makes sense to consider this operation for all skew left braces, although in
general it will not be associative.  The following problem appears
in~\cite{MR3818285}. 

\begin{problem}
	\label{pro:star_associative}
	Let $A$ be a left brace such that the operation $*$ is associative. Is $A$
	a two-sided brace?
\end{problem}

It is proved in~\cite[Proposition 2.2]{MR3818285} that the answer is positive
if the additive group of the left brace has no elements of order two. For skew
braces, Problem~\ref{pro:star_associative} can be answered in the negative,
see~\cite[\S1]{KSV}.  

A \emph{left ideal} of a skew left brace $A$ is a subgroup $L$ of the additive
group of $A$ such that $A*L\subseteq L$. An \emph{ideal} of $A$ is a left ideal
$I$ such that $a+I=I+a$ and $a\circ I=I\circ a$ for all $a\in A$. A skew left
brace $A$ is said to be \emph{simple} if it has no ideals different from
$\{0\}$ and~$A$. 

\begin{problem}
	\label{pro:simple}
	Classify finite simple skew left braces.
\end{problem}

Problem~\ref{pro:simple} is intensively studied for finite left braces, see for
example~\cite{MR3763276,BCJO18,MR3812099}. However, not much is known about
finite simple skew left braces that are not of abelian type.

\begin{problem}
	\label{pro:representations}
	Study representation theory of skew left braces.
\end{problem}

Rump defined modules over right braces in~\cite{MR2320986}. A different
interpretation of Problem~\ref{pro:representations} could be the following: To
study consequences of having skew left brace homomorphisms between two given
skew left braces.  

Left braces such that either the additive or the multiplicative group is
isomorphic to $\Z$ were classified in~\cite{CSV}. It is not clear how to extend
some of these results to skew left braces.  The operations $g^k + g^l =
g^{k+(-1)^kl}$ and $g^k\circ g^l = g^{k+l}$ turn the set $\{g^k:k\in\Z\}$ into
a skew left brace with multiplicative group isomorphic to $\Z$ and additive
group isomorphic to the infinite dihedral group.  It seems that there are no
other non-trivial skew left braces with multiplicative group isomorphic to
$\Z$.

\begin{problem}
	\label{pro:Z}
	Classify isomorphism classes of skew left braces with multiplicative group
	isomorphic to $\Z$.
\end{problem}

Problem~\ref{pro:Z} could be studied using techniques from factorizable groups,
see for example~\cite{MR1211633}. Similarly one could also ask for the
classification of skew left braces where one of its groups is isomorphic to the
infinite dihedral group.

Finite skew left braces with cyclic additive groups were classified by Rump
in~\cite{MR2298848,Rump}.  
In the same vein, it would be interesting to know the classification of skew
left braces with multiplicative group isomorphic to the quaternion group.  This
problem is not even solved for left braces.  Such left braces are called
\emph{quaternion left braces}. For $m\in\N$ let $q(4m)$ be the number of
isomorphism classes of quaternion left braces of size $4m$. 

\begin{problem}
	\label{pro:quaternion}
	Prove that 
	\[
		q(4m)=\begin{cases}
			2 & \text{if $m$ is odd,}\\
			7 & \text{if $m\equiv0\bmod8$,}\\
			9 & \text{if $m\equiv4\bmod8$,}\\
			6 & \text{if $m\equiv2\bmod8$ or $m\equiv6\bmod8$}.\\
		\end{cases}
	\]
\end{problem}

The conjectural formula for $q(4m)$ of Problem~\ref{pro:quaternion} was verified by
computer for all $m\leq128$. 

\begin{problem}
		Which finite abelian groups appear as the additive group of a
		quaternion left brace?
\end{problem}

For $m\in\{2,\dots,128\}$ the additive group of a quaternion left brace of order $4m$
is isomorphic to one of the following groups: 
\begin{equation*}
        C_{4m},\,
        C_{2m}\times C_2,\,
        C_m\times C_2\times C_2,\,
        C_m\times C_4,\,
        C_{m/2}\times C_2\times C_2\times C_2.
\end{equation*}
By inspection, $C_{m}\times C_2^2$ appears if 
$m\equiv2,4,6\bmod8$ and $C_m\times C_4$ and $C_{m/2}\times C_2^3$
appear if $m\equiv4\bmod8$.

\begin{problem}
    \label{con:Q}
	Is it true that there are seven classes of isomorphism of quaternion left braces of
	size $2^k$ for $k>4$? 
\end{problem}

It is interesting to study properties of groups appearing as the multiplicative
groups of skew left braces. Let us start with the case of left braces. A
finite group $G$ is said to be an \emph{involutive Yang--Baxter} group (or
simply a IYB group) if it is isomorphic to the multiplicative group of a finite
left brace.  In~\cite{MR1722951}, Etingof, Schedler and Soloviev proved that IYB
groups are always solvable.  A natural problem arises: Is every finite solvable
group the multiplicative group of a left brace?  The answer is negative, as it was
shown by Bachiller~\cite{MR3465351} following the ideas of
Rump~\cite{MR3291816}.  However, it would be interesting to find other
counterexamples.

\begin{problem}
	Which is the minimal cardinality of a solvable group that is not a IYB group?
\end{problem}

Some results related to IYB-groups can be found in~\cite{
MR3500774,
MR2584610,MR3177933,MR3447734,MR3116652}.  Related problems are the following:

\begin{problem}
	\label{pro:Cedo1}
	Is every nilpotent group of nilpotecy class two the multiplicative group
	of a left brace? 
\end{problem}

\begin{problem}
	\label{pro:Cedo2}
	Is every nilpotent group of nilpotecy class two the multiplicative group
	of a two-sided brace?
\end{problem}

\begin{problem}
	\label{pro:Cedo3}
	Which finite nilpotent groups are multiplicative groups of a
	two-sided (or left) brace?
\end{problem}

Problems~\ref{pro:Cedo1}--\ref{pro:Cedo3} appeared in the survey~\cite{cedo}.
Problem~\ref{pro:Cedo3} is interesting even in the particular case of groups of nilpotency class
$\leq3$. 

I learned the following related problem from Rump:

\begin{problem}
	\label{pro:IYBSylow}
	Is there an example of a non-IYB finite group where all the Sylow subgroups
	are IYB groups?
\end{problem}

The solution to the YBE one obtains from a skew left brace is always a
biquandle. Biquandles are non-associative structures useful in combinatorial
knot theory.  We refer to~\cite{MR3379534,MR2896084,MR3444659} for nice
introductions to the subject.

\begin{problem}
	Study knot invariants produced from skew left braces.
\end{problem}

As biquandles produce knot invariants and these invariants can be strengthened
by using quandle and biquandle homology, it is then natural to ask if the
homology of braces of~\cite{MR3530867} can be used in combinatorial knot
theory: 

\begin{problem}
	Is it possible to strengthen invariants produced from skew left braces by
	using brace homology?
\end{problem}

Left and right nilpotent left (and right) braces were defined by Rump in~\cite{MR2278047}.
Strongly nilpotent left braces were defined by Smoktunowicz in~\cite{MR3814340}.
These definitions extend to skew left braces, see~\cite{CSV}.  A skew brace
$A$ is said to be \emph{left nilpotent} if there exists $n\geq1$ such that
$A^n=0$, where $A^1=A$ and $A^{n+1}$ is defined as the subgroup $A*A^n$ of $(A,+)$
generated by $\{a*x:a\in A,x\in A^n\}$. Similarly $A$ is said to
be \emph{right nilpotent} if there exists $n\geq1$ such that $A^{(n)}=0$, where
$A^{(1)}=A$ and $A^{(n+1)}$ is defined as the subgroup $A^{(n)}*A$ of $(A,+)$ generated by
$\{x*a:x\in A^{(n)},a\in A\}$. 

The following problem of Smoktunowicz appears in~\cite{MR3765444} in an equivalent formulation:

\begin{problem}
	\label{pro:Agata}
	Let $G$ be a finite group which is the multiplicative group of some left
	brace. Is $G$ is the multipicative group of a right nilpotent left brace?
\end{problem}

Problem~\ref{pro:Agata} also makes sense for skew left braces of nilpotent type.

\begin{problem}
	\label{pro:CSV1}
	Are there simple two-sided skew braces of nilpotent type?
\end{problem}

A skew brace $A$ is said to be \emph{prime}
if for all non-zero ideals $I$ and $J$, the subgroup $I*J$ of $(A,+)$ generated by
$\{u*v:u\in I,v\in J\}$ is non-zero. 

\begin{problem}
	\label{pro:CSV2}
	Is every finite prime skew left brace of nilpotent type a simple skew left brace?
\end{problem}

\begin{problem}
	\label{pro:CSV3}
	Is every finite prime left brace a simple left brace?
\end{problem}

In ~\cite[\S5]{abundance} Ced\'o, Jespers and Okni\'nski found a
prime non-simple finite left brace. This example solves
Problems~\ref{pro:CSV2} and~\ref{pro:CSV3}.

\begin{problem}
	\label{pro:CSV4}
	Are there prime two-sided skew braces of nilpotent type?
\end{problem}

A skew left brace is said to be \emph{strongly nilpotent} if it is right and
left nilpotent.  A skew left brace $A$ is said to be \emph{strongly nil} if for
every element $a\in A$ there is a positive integer $n=n(a)$ such that any
$*$-product of $n$ copies of $a$ is zero.

\begin{problem}
	\label{pro:CSV5}
	Is every finite strongly nil skew left brace a strongly nilpotent skew left brace?
\end{problem}

For a skew left brace $A$ let $\rho_1(a)=a$ and $\rho_{k+1}(a)=\rho_k(a)*a$ for
$n\geq1$.  The skew left brace $A$ is said to be \emph{right nil} if there
exists a positive integer $n$ such that $\rho_n(a)=0$ for all $a\in A$. 

\begin{problem}
	\label{pro:CSV6}
	Is every finite right nil skew left brace a right nilpotent skew left brace?
\end{problem}

Radical rings are the source of several other problems for skew left braces. 
For example, it might make sense
to discuss an analog of the K\"othe conjecture in the context of skew left braces.
For rings one of the formulations of the conjecture is the following: the sum of two nil left ideals in a
ring is a left nil ideal. The conjecture was formulated arround 1930 and it is
still open, see~\cite{MR3439122,MR1879880}. 

\begin{problem}
	\label{pro:CSV7}
	Is there an analog of the K\"othe conjecture for skew left braces?
\end{problem}

Problems~\ref{pro:CSV1},~\ref{pro:CSV2},~\ref{pro:CSV3},~\ref{pro:CSV4},~\ref{pro:CSV5}
and~\ref{pro:CSV6} appear in~\cite{CSV} and~\cite{KSV}. K\"othe conjecture has
been shown to be true for various classes of rings such as polynomial identity
rings, right Noetherian rings and radical rings. 

Bachiller observed the connection between skew braces and Hopf--Galois
extensions. This connection is precisely described in the appendix of~\cite{MR3763907}.
In the theory of Hopf--Galois extensions the following problems are natural,
see~\cite{MR3425626}:

\begin{problem}
	\label{pro:SS}
	Is there a skew left brace with solvable additive group but non-solvable
	multiplicative group?
\end{problem}

\begin{problem}
	\label{pro:NS}
	Is there a skew left brace with non-solvable additive group but nilpotent
	multiplicative group?
\end{problem}

Problems~\ref{pro:SS} and~\ref{pro:NS} also appeared in~\cite{MR3763907}
and~\cite[Problem 19.90]{kourovka}. In general, it is interesting to know
relations between the multiplicative and the additive group of a skew left
brace. Several results in this direction can be found in the theory of
Hopf--Galois extensions, see for example~\cite{MR2011974,MR3030514,MR3425626}.
The following problem is mentioned in~\cite{LV3} and~\cite[Problem
19.49]{kourovka}:

\begin{problem}
	\label{pro:orderable}
	Let $A$ be a skew left brace with left-orderable multiplicative group. Is the
	additive group of $A$ left-orderable?
\end{problem}

%
Timur Nasybullov constructed an example of a left brace that answers
Problem~\ref{pro:orderable} negatively (private communication). Another example that answers Problem~\ref{pro:orderable} appears in~\cite{CSV}.

In~\cite{CSV} it is proved that skew left braces generated (as a skew left
brace) by a single element yield indecomposable solutions of the YBE. It is
then natural to ask if all finite indecomposable solutions are of this type:

\begin{problem}
	\label{pro:AA}
	Is it true that for each indecomposable involutive solution $(X,r)$
	there exists a one-generator left brace $A$ generated by $x$ such that $X$ is
	of the form $X=\{ax+x:a\in A\}$?
\end{problem}

Problem~\ref{pro:AA} appears in~\cite{MR3771874} for left braces. 

Rump proved that braces are equivalent to linear cycle sets.  A \emph{linear
cycle set} is a triple $(A,+,\cdot)$, where $(A,+)$ is an abelian group 
and $(A,\cdot)$ is a cycle set 
such that
$a\cdot (b+c)=(a\cdot b)+(a\cdot c)$ and 
$(a+b)\cdot c=(a\cdot b)\cdot (a\cdot c)$ 
hold for all $a,b,c\in A$.

A theory of dynamical extensions of cycle sets was used in~\cite{MR3437282} to
produce a counterexample to Gateva--Ivanova conjecture.  

\begin{problem}
	\label{pro:dynamical}
	Develop the theory of dynamical extensions of linear cycle sets.
\end{problem}

To study non-involutive solutions one replaces cycle sets by skew cycle sets.
According to~\cite[\S6.1]{MR3763907}, a skew cycle set is defined as a linear
cycle set but where the commutativity of the group $(A,+)$ is not assumed.
Problem~\ref{pro:dynamical} can be stated for skew braces and skew 
cycle sets.


Almost all of the questions for skew left braces make sense in the particular
and highly interesting case of $k$-linear left braces, which are left braces
where the additive group is a $k$-vector space compatible with the
multiplicative group.  Those braces were introduced by Catino and Rizzo
in~\cite{MR2486886} as \emph{circle algebras}.

There are also several questions on one-generator skew left braces, as those
skew left braces are maybe easier to study than skew left braces. In
particular, one could also ask most of the previous questions on skew left braces
for one-generator skew left braces.

\subsection*{Acknowledgements}

This work is supported by PICT-201-0147 and MATH-AmSud 17MATH-01. The author
thanks F. Ced\'o, T. Gateva--Ivanova, E. Jespers, V. Lebed, K. Nejabati Zenouz,
J. Okni\'nski, A.  Smoktunowicz, W.  Rump, P. Zadunaisky for discussions,
comments and problems.

\nocite{MR3478858}

\bibliographystyle{abbrv}
\bibliography{refs}

\end{document}